\documentclass[conference,letterpaper,10pt,twocolumn,tbtags]{IEEEtran}

\usepackage{comment}

\usepackage{ifpdf}
\usepackage[
  pdfstartview={FitH},
	pdfpagemode={None},
  pdfauthor={Aditya Mahajan and Demosthenis Teneketzis},
  pdftitle={Optimal Performance of Feedback Control Systems with Limited
  Communication over Noisy Channel},
  pdfsubject={IEEE 45th Conference on Decison and Control},
  pdfkeywords={}
 ]{hyperref}

\usepackage{cite}
\makeatletter
\def\citedash{\hbox{--}\penalty\@m}
\makeatother

\usepackage[T1]{fontenc}
\usepackage[latin1]{inputenc}

\usepackage{graphicx}
\usepackage{emp}
\renewcommand\empTeX{\documentclass[onecolumn]{IEEEtran}}

\ifpdf
\else
\fi

\usepackage[cmex10]{empheq}
\usepackage{nccmath,braket,amssymb}
\usepackage{maybemath}
\usepackage{dsfont,mathrsfs}
\let\qedhere\relax

\newcounter{mytempeqncnt}

\newtheorem{problem}{Problem}
\newtheorem{theorem}{Theorem}
\newtheorem{definition}{Definition}
\newtheorem{lemma}{Lemma}

\DeclarePairedDelimiter\abs{\lvert}{\rvert}

\newcommand\E[1]{\mathds E{\Set{#1}}}

{\catcode`\|=\active \gdef\conditional#1{\begingroup 
\mathcode`\|32768\let|\BraVert\left({#1}\right)\endgroup}
}

\let\PPr\Pr
\renewcommand\Pr[2][]{\PPr{#1}{\conditional{#2}}}
\makeatother

\newcommand*\defined{\triangleq}

\newcommand*\I[1]{\;\mathds 1\left[\,#1\,\right]}

\newcommand*\sC{\mathscr C}
\newcommand*\sF{\mathscr F}

\newcommand*\sG{\mathscr G}
\newcommand*\sL{\mathscr L}
\renewcommand*\P[1]{\mathcal P^{#1}}

\newcommand*\J{\mathcal J}

\newcommand*\M{\mathcal M}
\newcommand*\N{\mathcal N}

\renewcommand*\S{\mathcal S}
\newcommand*\U{\mathcal U}

\newcommand*\W{\mathcal W}
\newcommand*\X{\mathcal X}
\newcommand*\Y{\mathcal Y}
\newcommand*\Z{\mathcal Z}

\newcommand*\RR{\mathds R}

\newcommand\PXM[2]{\ensuremath{P_{X_{#1},M_{#2}}}}
\newcommand\tPXM[2]{\ensuremath{\widetilde P_{X_{#1},M_{#2}}}}

\newcommand\PXMp[2]{\ensuremath{P^+_{X_{#1},M_{#2}}}}
\newcommand\tPXMp[2]{\ensuremath{\widetilde P^+_{X_{#1},M_{#2}}}}

\newcommand\PXSM[3]{\ensuremath{P_{X_{#1},S_{#2},M_{#3}}}}

\usepackage{xspace}
\newcommand*\mbmo{ member by member optimal (m.b.m.o.)%
\gdef\mbmo{ m.b.m.o.\@\xspace}\xspace}

\newcommand*\pomdp{partially observed Markov decision problems (\textsc{pomdp})\xspace%
\gdef\pomdp{\textsc{pomdp}\xspace}}

\newcommand*\POsystem{ \ensuremath{(\X,\allowbreak \W,\allowbreak \M, 
\allowbreak \Z,\allowbreak \N,\allowbreak \Y,\allowbreak \U,\allowbreak  
P_{X_1},\allowbreak  P_W,\allowbreak  P_N,\allowbreak  f,\allowbreak  
h,\allowbreak  \rho,\allowbreak  T)}}

\newcommand*\IOsystem{ \ensuremath{(\X,\allowbreak \W,\allowbreak \widehat 
\N,\allowbreak \S,\allowbreak \M, \allowbreak \Z,\allowbreak \N,\allowbreak 
\Y,\allowbreak \U,\allowbreak  P_{X_1},\allowbreak  P_W,\allowbreak 
P_{\widehat N}, \allowbreak  P_N,\allowbreak  f,\allowbreak \mu, \allowbreak  
h,\allowbreak  \rho,\allowbreak  T)}}

\newcommand*\design[1][]{\ensuremath{c^{#1},l^{#1},g^{#1}}}

\def\today{\ifcase\month\or
 January\or February\or March\or April\or May\or June\or
 July\or August\or September\or October\or November\or December\fi
 \space\oldstylenums{\number\day}, \oldstylenums{\number\year}}

\newif\ifProcessMP
\ProcessMPfalse

\newcommand\startfloatequation[1]%
  {\begin{figure*}[!t]
    \normalsize
    \setcounter{mytempeqncnt}{\value{equation}}
    \setcounter{equation}{#1}}

\newcommand\stopfloatequation[1]
  { 
    \setcounter{equation}{\value{mytempeqncnt}}
    \hrulefill
    \vspace*{4pt}
  \end{figure*}
  \addtocounter{equation}{#1}}

\begin{document}

\title{Optimal Performance of Feedback Control Systems with Limited 
Communication over Noisy Channels} 
\author{\authorblockN{Aditya Mahajan and Demosthenis 
Teneketzis~\IEEEmembership{Fellow,~IEEE}}%
\authorblockA{Department of Electrical Engineering and Computer Science\\
University of Michigan, Ann Arbor, MI 48109--2122, USA\\
Email: \{adityam,teneket\}@eecs.umich.edu}
}
\maketitle

\begin{abstract}
  A discrete time stochastic feedback control system with a noisy
  communication channel between the sensor and the controller is considered.
  The sensor has limited memory. At each time, the sensor transmits encoded
  symbol over the channel and updates its memory. The controller receives a
  noisy version of the transmitted symbol, and generates a control action
  based on all its past observations and actions. This control action action
  is fed back into the system. At each stage the system incurs an
  instantaneous cost depending on the state of the plant and the control
  action.  The objective is to choose encoding, memory updating and control
  strategies  to minimize the expected total costs over a finite horizon, or
  the expected discounted cost over an infinite horizon, or the expected
  average cost per unit time over an infinite horizon.  For each case we
  obtain a  sequential decomposition of the optimization problem. The results
  are extended to the case when the sensor makes an imperfect observation of
  the state of the system.  
\end{abstract}

\section{Introduction}

\noindent Recent advances in network and communication technologies have led
to an increasing interest in networked control systems (\textsc{ncs})~(see the
papers in~\cite{NCS:T-AC}), in particular, the limitations imposed upon
feedback control by the presence of a communication channel in the loop. Most
researchers have concentrated on stability analysis of the system.  The
problem of stabilization of a plant with finite data rate feedback  was
investigated in~\cite{Delchamps:1990, WongBrockett:1999, Baillieul:2001,
Baillieul:2002, EliaMitter:2001, BrockettLiberzon:2000, PetersonSavkin:2001,
IshiiFrancis:2003, Liberzon:2003, NairEvans:2000, NairEvans:2003,
NairEvans:2004, NairEvansMareelsMoran:2004, MartinsDahlehElia:2004}.
\textsc{lqg} stability of deterministic and stochastic systems under various
communication constraints (rate limited channels, noisy channels with input
power constraint, etc.) was considered
in~\cite{BorkarMitterTatikonda:2001, Tatikonda:phd, Tatikonda:2004,
Tatikonda:2004a, Tatikonda:2004b, MatveevSavkin:2004, Braslavsky:2005}.  Performance limitation in terms of
lower bounds on the separation of differential entropy rates was investigated
in~\cite{Martins1, Martins2}. However certain applications require performance
metrics more general than asymptotic metrics of stability and separation
of differential entropy. In this paper we consider the class of additive
performance metrics, where the total cost is the sum of costs along the entire
path.

In problems with asymptotic performance metrics, transient behavior need not
be optimal, thus strong performance bounds can be derived by using asymptotic
results from probability theory, information theory and classical control
theory.  However, in problems with more general performance metrics, transient
behavior needs to be optimal. To the best of our knowledge, performance
analysis of such problems has not been addressed in the literature. We
identify algorithms to obtain optimal strategies; but we have not been able to
find expressions for optimal performance or performance bounds.

We consider a discrete-time feedback control system with a communication
channel between the sensor and the controller, shown in \autoref{fig:system}.
Such  problems arise when the plant and the controller are geographically
separated.  We are interested in problems in which the sensor has limited
resources while the controller has no resource constraint; we model the sensor
as an encoder with finite memory, the channel between the sensor and the
controller as noisy, and the channel between the controller and the system as
noiseless\footnote{In the sequel we show that this assumption does not entail
any loss of generality.}.  At each stage the system incurs an instantaneous
cost depending on the state of the plant and the control action. The objective
is to choose encoding, memory updating and control strategies to minimize the
expected total cost over a finite horizon, or expected discounted cost over an
infinite horizon, or expected average cost per unit time over an infinite
horizon.

The key contribution of this paper is providing a methodology for determining
jointly optimal real-time encoding, memory updating and control strategies for
feedback control systems with limited communication over noisy channels. The
methodology applies to  general non-linear stochastic systems,
with an arbitrary additive performance criteria. 

The remainder of this paper is organized as follows. We formulate the
performance analysis of feedback control systems with limited communication
over noisy channels as a decentralized stochastic optimization problem. To
illustrate the key concepts associated with our solution methodology we first
consider in \autoref{sec:finite} the finite horizon problem, for which we
establish structural results of optimal controller and present a methodology for joint
optimization of the encoding, memory updating and control strategies. In
\autoref{sec:infinite} we extend the methodology to infinite horizon problems.
We discuss computational issues for obtaining numerical solution of the
dynamic programming algorithms of finite and infinite horizon problems in
\autoref{sec:computation}.  The feedback control problem when the encoder has
imperfect observation of the state of the plant is considered in
\autoref{sec:partial}. We conclude in \autoref{sec:conclusion}.

\paragraph*{Notation}
We use uppercase letters $(X,Y,Z)$ to denote random variables and lowercase 
letters to denote their realizations $(x,y,z)$.  When we represent a function 
of random variables as a random variable $(\PXM tt, \PXM {t+1}t, \PXMp 
{t+1}t)$, a tilde above the variable denotes its realization $(\tPXM tt, 
\tPXM {t+1}t, \tPXMp {t+1}t)$.  When we use Greek letters to represent a 
random variable $(\pi^c, \pi^l, \pi^g)$, a tilde above the variable denotes 
its realization $(\tilde \pi^c, \tilde \pi^l, \tilde \pi^g)$. We also use the 
short hand notation of $x^t$ to present the sequence $x_1,\dots,x_t$ and 
similar notation for random variables and functions.

\section{The Finite Horizon Problem}\label{sec:finite}

\subsection{Problem Formulation}\label{sec:prob}
\begin{figure}[tb!h]
  \centering
  \resizebox{0.95\linewidth}{!}{
  \empuse{feedback}}
    \caption{Feedback control system with noisy communication}
  \label{fig:system}
\end{figure}
Consider a discrete time feedback control system as shown in 
\autoref{fig:system} which operates for a horizon $T$. The state evolution is 
given by
\begin{equation}
  X_{t+1} = f(X_t, U_t, W_t),
  \label{eq:state}
\end{equation}
where $f$ is the \emph{system evolution function}. The variables $X_t, U_t, 
W_t$ denote the state of the system, the control action and the plant disturbance 
respectively, at time $t$. We assume that all variables are discrete. For all 
$t$, $X_t$ takes values in $\X \defined \{ 1, 2, \dots, \abs\X \}$, $U_t$ 
takes values in $\U \defined \{1, 2, \dots, \abs\U\}$ and $W_t$ takes values 
in $\W \defined \{1, 2, \dots, \abs\W\}$.  The initial state $X_1$ is a 
random variable with PMF $P_{X_1}$. The random variables $W_1, \dots, W_T$ 
are i.i.d.\@ with PMF $P_W$ and are also independent of $X_1$.

The sensor, consisting of an encoder and a memory, makes perfect observations 
of the state of the system. At each time instant $t$, the encoder generates an
encoded symbol $Z_t$, taking values in $\Z \defined \{1, \dots, \abs\Z\}$, as
follows
\begin{equation}
  Z_t = c_t(X_t,M_{t-1}),
  \label{eq:encoder}
\end{equation}
where $c_t$ is the \emph{encoding function} at time $t$ and $M_{t-1}$ denotes 
the content of the sensor's memory at $t-1$. $M_t$ takes values in $\M 
\defined \{1, \dots, \abs\M\}$ and is updated according to
\begin{equation}
  M_t = l_t(X_t, M_{t-1}),
  \label{eq:memory}
\end{equation}
where $l_t$ is the \emph{memory update function} at time $t$. Observe that the 
sensor has a finite size memory and though it makes perfect observations of 
the state of the system, it can not store all the past observations. Thus, it 
does not have perfect recall and at each stage it must selectively shed 
information.

The encoded symbol $Z_t$ is transmitted over a noisy communication channel and 
a channel output $Y_t$ is generated according to
\begin{equation}
  Y_t = h(Z_t, N_t),
  \label{eq:channel}
\end{equation}
where $h$ is the channel and $N_t$ denotes the channel noise. $Y_t$ takes 
values in $\Y \defined \{1,\dots, \abs\Y \}$ and $N_t$ takes values in $\N 
\defined \{1,\dots, \abs\N\}$. The sequence of random variables 
$N_1,\dots,N_T$ is i.i.d.\ with given PMF $P_N$. $N_1,\dots,N_T$ are also 
independent of $X_1,W_1,\dots,W_T$.

The controller observes the channel outputs and generates a control action 
$U_t$ as follows
\begin{equation}
  U_t = g_t(Y^t,U^{t-1}),
  \label{eq:controller}
\end{equation}
where $g_t$ is the \emph{control law} at time $t$. $U_t$ takes values in $\U 
\defined \{1,\dots, \abs\U\}$. A uniformly bounded cost function $\rho: 
\X\times\U \to [0,K]$, $K < \infty$ is given. At each $t$, an instantaneous 
cost $\rho(X_t,U_t)$ is incurred.

The collection $\POsystem$ is called a \emph{perfect observation system}. The 
choice of (\design), $c \defined (c_1,\dots,c_T)$, $l \defined (l_1, \dots, 
l_T)$, $g \defined (g_1, \dots, g_T)$, is called a \emph{design}.

The performance of a design, quantified by the expected total cost under 
that design, is given by
\begin{equation}
  \J_T(\design) \defined \E{\sum_{t=1}^{T} \rho(X_t, U_t) | \design },
  \label{eq:cost}
\end{equation}
where the expectation in \eqref{eq:cost} is with respect to a joint measure on 
$(X_1,\dots, X_T, U_1, \dots, U_T)$ generated by $P_W, P_N, f, h$ and the 
choice of design (\design). We are interested in the following optimization 
problem:
\begin{problem}\label{prob:basic}
  Given a perfect observation system $\POsystem$, choose a design 
  $(\design[*])$ such that
  \begin{equation}
    \J_T(\design[*]) = \J^*_T \defined \min_{\design \in \sC^T \times \sL^T 
    \times \sG^T } \J_T(\design),
  \end{equation}
  where $\sC^T \defined \sC \times \dots \times \sC$ ($T$ times), $\sC$ is the 
  space of functions from $\X\times\M$ to $\Z$, $\sL^T \defined \sL \times 
  \dots \times \sL$ ($T$ times), $\sL$ is the space of functions from 
  $\X\times\M$ to $\M$, $\sG^T \defined \sG_1 \times \dots \times \sG_T$, and 
  $\sG_t$ is the space of functions from $\Y^t\times\U^{t-1}$ to $\U$.
\end{problem}

\paragraph{Remark}
There is no loss of generality in assuming a noiseless feedback channel.
Suppose there is noise in the feedback channel, and the input to the system is
$\widehat U_t$ is a noisy version of $U_t$ given by
\begin{equation}
  \widehat U_t = \hat h(U_t, \widehat N_t)
\end{equation}
where $\hat h$ is the feedback channel and $\widehat N_t$ denotes the noise in
the feedback channel. $\widehat N_1, \dots, \widehat N_T$ is a sequence of
independent variables that is also independent of $X_1, W_1\dots, W_T$ and
$N_1, \dots, N_T$. This model can be transformed into one equivalent to
\eqref{eq:state}--\eqref{eq:controller} by setting
\begin{align}
  \widehat W_t &= (W_t, \widehat N_t), \\
  X_{t+1} &= f\left( X_t, h(U_t, \widehat N_t), W_t \right) \defined \hat f 
  (X_t, U_t, \widehat W_t).
\end{align}
Thus, without loss of generality we can assume a noiseless feedback channel.

\subsection{Salient Features of the Problem}\label{sec:salient}

\autoref{prob:basic} is a decentralized multi-agent stochastic optimization
problem. The agents---the sensor and the controller---share a common objective
of minimizing the expected total cost. They have access to different (and
non-nested) information about underlying state of nature. Furthermore, the
actions taken by an agent at any instant of time affects the observations of
the other agent at future time instants. Thus the problem is a sequential (in
the sense of~\cite{Witsenhausen:1971a}) \emph{dynamic team} with strictly
non-classical information structure~\cite{Witsenhausen:1971}.  
Dynamic teams are, in general, functional optimization problems having a
complex interdependence among the decision rules~\cite{Ho:1980}.  This
interdependence leads to non-convex (in policy space) optimization problems
that  are hard to solve (see~\cite{Witsenhausen:1968} for an example). 
Identifying an information state sufficient for performance
evaluation~\cite{Witsenhausen:1976,KumarVaraiya:1986} is a
key step in obtaining a sequential decomposition of such problems.
To obtain a sequential decomposition of \autoref{prob:basic}, we proceed as
follows. First, we derive structural properties of optimal controllers. Using
these structural results we transform \autoref{prob:basic} into an equivalent
optimization problem and identify information states for this equivalent
problem. This yields a sequential decomposition for \autoref{prob:basic} along
with a dynamic programming algorithm to obtain an optimal design.

\subsection{Structural Results}
In this section we present structural properties of optimal controllers.  For 
this purpose we define the following.
\begin{definition}\label{def:beliefs}
  Let \PXM tt, \PXM {t+1}t and \PXMp {t+1}t be random vectors defined as 
  follows:
  \begin{align}
      &\PXM tt(x,m)  \nonumber\\
        =\Pr{ X_t = x, M_t = m | Y^t, U^{t-1}, c^t, l^t, g^{t-1}},
        \span \omit \\
    &\PXM {t+1}t(x,m)  \nonumber\\
      = \Pr{ X_{t+1} = x, M_t = m | Y^t, U^t, c^t, l^t, g^t},
        \span \omit \\
    &\PXMp{t+1}t(x,m)  \nonumber\\ 
      \quad = \Pr{ X_{t+1} = x, M_t = m | Y^{t+1}, U^t, c^{t+1}, l^t, g^t}.
      \span\omit
    \end{align}
\end{definition}
For any particular realization $y^t, u^{t-1}$ and arbitrary (but fixed) 
choice of $c^t$, $l^{t}$ and $g^{t-1}$, the realization of \PXM tt, denoted 
by \tPXM tt, is a PMF on $(X_t, M_t)$. If $(Y^t, U^{t-1})$ is a random vector 
and $c^t, l^{t}, g^{t-1}$ are arbitrary (but fixed) functions, then \PXM tt 
is a random vector belonging to $\P{\X\times\M}$, the space of PMFs on 
$\X\times\M$.  Similar interpretations hold for \PXM {t+1}t and \PXM {t+1}t.


These beliefs given by \autoref{def:beliefs} are related as
follows:
\begin{lemma}\label{lemma:beliefs}
  For each stage $t$, there exists a deterministic functions $\psi$, $\phi$,
  and $\nu$ such that 
  \begin{align}
    \PXM {t+1}t &= \psi( \PXM tt, U_t ), 
    \label{eq:psi} \\
    \PXMp {t+1}t &= \phi( \PXM {t+1}t, c_{t+1} ),
    \label{eq:phi} \\
    \PXM {t+1}{t+1} &= \nu( \PXMp {t+1}t, l_{t+1} ).
    \label{eq:nu}
  \end{align}
\end{lemma}
\begin{proof}
  Consider a component of \tPXM {t+1}t, given by~\eqref{eq:A_a} at the top of
  the page.
  \startfloatequation{13}
     \begin{gather}
         \tPXM {t+1}t(x,m)
         = \frac{ \Pr{X_{t+1} = x, M_t = m, U_t = u_t | y^t, u^{t-1}, c^t, 
         l^t, g^t }}
         { \sum\limits_{x',m' \in \X \times \M}
           \Pr{X_{t+1} = x', M_t = m', U_t = u_t | y^t, 
           u^{t-1}, c^t, l^t, g^t }}.
       \label{eq:A_a}
     \end{gather}
    \stopfloatequation{1}

        Now consider
      \begin{equation}
        \begin{split}
          &\Pr{X_{t+1} = x, M_t = m, U_t = u_t | y^t, u^{t-1}, c^t, l^t, g^t 
          }  \\
          &= \sum_{x_t \in \X} \Pr{X_t = x_t, M_t = m | y^t, u^{t-1}, c^t, l^t, 
          g^t} \\
          &\quad\times
          \Pr{U_t = u_t | y^t, u^{t-1}, c^t, l^t, g^t, x_t, M_t = m} \\
          &\quad \times \Pr{X_{t+1} = x | x_t, M_t = m, y^t, u^t, c^t, l^t, 
          g^t} \\
          &\stackrel{(a)}= \Pr{u_t | y^t, l^{t-1}, g_t} 
          \\
          &\times
          \sum_{x_t \in \X}
          \Big[
          \Pr{X_t = x_t, M_t = m | y^t, u^{t-1}, c^t, l^t, g^{t-1}} \\
          \times
          \Pr{X_{t+1} = x | x_t, u_t}\Big] \span\omit\\
          &=\Pr{u_t | y^t, l^{t-1}, g_t } \\
          &\quad\times\sum_{x_t \in \X}
          \Big[ \tPXM tt (x_t,m) \Pr{X_{t+1} = x | x_t, u_t} \Big],
          \end{split}
        \label{eq:A_b}
      \end{equation}
      where $(a)$ follows from~\eqref{eq:state} and~\eqref{eq:controller}.  
      Combine~\eqref{eq:A_a} and~\eqref{eq:A_b} and cancel 
      $\Pr{u_t|y^t,l^{t-1},g_t}$ from the numerator and the 
      denominator of~\eqref{eq:A_a}, giving
      \begin{equation}
        \tPXM {t+1}t  = \psi(\tPXM tt, u_t),
      \end{equation}
      where $\psi$ is given by~\eqref{eq:A_a} and~\eqref{eq:A_b}.
    
      Consider a component of \tPXMp {t+1}t, given by~\eqref{eq:b_A} on the
      top of the page.
      \startfloatequation{16}
        \begin{gather}
             \tPXMp {t+1}t(x,m) 
            = \frac{ \Pr{X_{t+1} = x, M_t = m, Y_{t+1} = y_{t+1} | y^t, u^t, 
            c^{t+1}, l^t, g^t}}
            { \displaystyle \sum_{x',m' \in \X \times \M} \Pr{X_{t+1} = x', M_t = m', Y_{t+1} = 
            y_{t+1} | y^t, u^t, c^{t+1}, l^t, g^t}}.
          \label{eq:b_A}
        \end{gather}
      \stopfloatequation{1}

      Now consider,
      \begin{equation}
        \begin{split}
          &\mathrlap{
          \Pr{X_{t+1} = x, M_t = m, Y_{t+1} = y_{t+1} | y^t, u^t, c^{t+1}, 
          l^t, g^t}} \\
          &= \Pr{X_{t+1} = x, M_t = m | y^t, u^t, c^{t+1}, l^t, g^t} \\
          & \quad\times
          \PPr\bigl( Y_{t+1} = y_{t+1} \bigm| y^t, u^t, c^{t+1}, \\
           l^t, g^t, X_{t+1} = x, M_t = m \bigr) \span \omit \\
          & \stackrel{(b)}=\Pr{X_{t+1} = x, M_t = m | y^t, u^t, c^t, l^t, 
          g^t}\\
          & \quad \times
          \Pr{Y_{t+1} = y_{t+1} | X_{t+1} = x, M_t = m , c_{t+1}} \\
          &= \tPXM {t+1}t (x,m) \\
          &\quad\times \Pr{Y_{t+1} = y_{t+1} | X_{t+1} = x,
          M_t = m, c_{t+1}},
        \end{split}
        \label{eq:b_B}
      \end{equation}
      where $(b)$ follows from~\eqref{eq:encoder} and~\eqref{eq:channel}.  
      Combining~\eqref{eq:b_A} and~\eqref{eq:b_B} we have
      \begin{equation}
        \tPXMp {t+1}t = \phi(\tPXM {t+1}t, c_{t+1}),
      \end{equation}
      where $\phi(\cdot)$ is given by~\eqref{eq:b_A} and~\eqref{eq:b_B}.

    Consider a component of \tPXM {t+1}{t+1},
      \begin{equation}
        \begin{split}
          & \tPXM {t+1}{t+1}(x,m) \\
          & =\Pr{X_{t+1} = x, M_{t+1} = m | y^{t+1}, u^t, c^{t+1}, 
           l^{t+1},g^t} \\
          &= \sum_{m_t \in \M}
          \mathrlap{
          \Pr{X_{t+1} = x, M_t = m_t | y^{t+1}, u^t, 
          c^{t+1}, l^{t+1},g^t}} \\
          & \quad \times\PPr\bigl(M_{t+1} = m \bigm| X_{t+1} = x, M_t = m_t, y^{t+1}, u^t, \\
          c^{t+1}, l^{t+1},g^t\bigr) \span\omit \\
          &\stackrel{(c)}= \sum_{m_t \in \M} 
          \mathrlap{
          \Pr{X_{t+1} = x, M_t = m_t | y^{t+1}, u^t, c^{t+1}, l^t,g^t}}\\
          &\quad  \times \Pr{M_{t+1} = m | X_{t+1} = x, M_t = m_t, l_{t+1}} \\
          &\quad  \times \sum_{m_t \in \M} \tPXM {t+1}t (x,m_t) 
          \I{m = l_{t+1}(x,m_t)} \\
          &\defined \nu (\tPXM {t+1}t, l_{t+1}),
        \end{split}
      \end{equation}
      where $(c)$ follows from~\eqref{eq:memory} and $\I{\cdot}$ is the
      indicator function.
\end{proof}
The above relationships between the controller's beliefs lead to the 
structural results of the optimal controllers.

\begin{theorem}
  \label{thm:structural}
  Consider \autoref{prob:basic} for any arbitrary (but fixed) encoding and 
  memory update strategies $c$ and $l$, respectively. Then, without loss of 
  optimality, we can restrict attention to control laws of the form
  \begin{equation}
    U_t = g_t(\PXM tt).
    \label{eq:structural}
  \end{equation}
\end{theorem}
\begin{proof}
  Equations \eqref{eq:psi}--\eqref{eq:nu} of \autoref{lemma:beliefs} can be 
  combined to obtain
  \begin{equation}
    \begin{split}
      \PXM {t+1}{t+1} &= \nu \left( \phi\big(\psi(\PXM tt, U_t), c_{t+1} 
      \big), l_{t+1} \right) \\
      & \defined \mu(\PXM tt, U_t, c_{t+1}, l_{t+1}).
    \end{split}
  \end{equation}
  Thus for any fixed $c$ and $l$, \PXM tt is a controlled Markov process with 
  control action $U_t$. Further, the expected instantaneous cost can be 
  written as
  \begin{equation}\label{eq:equiv_cost}
    \begin{split}
      \hspace{1em} & \hspace{-1em}
    \E{\rho(X_t, U_t) | y^t, u^t, c^t, l^t, g^t} \\
    &= \sum_{x_t,m_t \in \X
    \times \M} \rho(x_t, 
    u_t) \PXM tt (x_t,m_t) \\
    &\defined \hat \rho(\tPXM tt, u_t).
  \end{split}
  \end{equation}
  There is a subtle technicality in the first step of~\eqref{eq:equiv_cost}.
  See~\cite{MahajanTeneketzis:2006b} for details. Hence we have a perfectly
  observed stochastic process $\{\PXM tt,\ t=1,\dots,T\}$ with control action
  $U_t$ and instantaneous cost $\hat \rho(\PXM tt, U_t)$. From Markov decision
  theory~\cite{KumarVaraiya:1986} we know that there is no loss of optimality
  in restricting attention to control laws of the form~\eqref{eq:structural}.
\end{proof}

\subsubsection{Implication of the structural results}
\autoref{thm:structural} implies that at each stage $t$, without loss of 
optimality, we can restrict attention to controllers belonging to 
the family $\sG_s$ of functions from $\P{\X\times\M}$ to $\U$. Thus at each stage we 
can optimize over a fixed (rather than time-varying) domain. Thus 
\autoref{prob:basic} is equivalent to the following problem:
\begin{problem}
  \label{prob:structural}
  Given a perfect observation system $\POsystem$, choose a design 
  (\design[*]) that is optimal with respect to the performance criterion of 
  \eqref{eq:cost}, i.e.,
  \begin{equation}
    \J_T(\design[*]) = \J^*_T \defined \min_{\design \in \sC^T \times \sL^T 
    \times \sG_s^T } \J_T(\design),
  \end{equation}
  where $\sG_s^T \defined \sG_s \times \dots \times \sG_s$ ($T$ times).
\end{problem}
Thus we have an optimization problem in which the action space is not changing 
with time. In the next section we provide a sequential decomposition of 
\autoref{prob:structural}.

\subsection{Joint Optimization}\label{sec:joint}
In this section, we identify information states sufficient for performance 
evaluation of \autoref{prob:structural}, resulting in its sequential 
decomposition. \autoref{prob:structural} is equivalent to 
\autoref{prob:basic}, hence we also obtain a sequential decomposition of  
\autoref{prob:basic}. The intuition behind our approach is as follows. As 
mentioned in \autoref{sec:salient}, the agents act in a sequential manner.  Let 
$\pi^c_t, \pi^l_t, \pi^g_t$ be the information states of the encoder, memory 
update and controller respectively. For these to be valid states, they must 
satisfy the property
\begin{equation*}
  \cdots \rightarrow \pi^c_t \xrightarrow{c_t} \pi^l_t \xrightarrow{l_t} 
  \pi^g_t \xrightarrow{g_t} \pi^c_{t+1} \rightarrow \cdots,
\end{equation*}
that is, at each time instant $t$, $\pi^l_t$ can be determined from $\pi^c_t$ and 
$c_t$, $\pi^g_t$ can be determined from $\pi^l_t$ and $l_t$, and 
$\pi^c_{t+1}$ can be determined from $\pi^g_t$ and $g_t$. This ensures that 
$\pi^c_t, \pi^l_t, \pi^g_t$ are \emph{information states} in the sense 
of~\cite{KumarVaraiya:1986}.  However, a system can have more than one 
information state, and  not all of them are sufficient for performance 
evaluation (see~\cite{Witsenhausen:1976}). To be sufficient for performance 
evaluation, the information states must \emph{absorb/summarize} the effect of 
past decision rules on the expected future cost, that is,
\begin{equation}
  \begin{split}
    \hspace{1em}&\hspace{-1em}
    \E{\sum_{s=t}^{T}\rho(X_s,U_s) | \design} \\
    &=
    \E{\sum_{s=t}^{T}\rho(X_s, U_s) | \pi^c_t, c_t^T, l_t^t, g_t^T } \\
    &=
    \E{\sum_{s=1}^{T}\rho(X_s, U_s) | \pi^l_t, c_{t+1}^T, l_t^t, g_t^T } \\
    &=
    \E{\sum_{s=t}^{T}\rho(X_s, U_s) | \pi^g_t, c_{t+1}^T, l_{t+1}^t, g_t^T }
  \end{split}
\end{equation}
Furthermore, to extend the results of the finite horizon problem to infinite 
horizon problems, we want the domain of information states to be 
time-invariant.

The following information states satisfy the above requirements.
\begin{definition}
  \label{def:info_states}
  Let $\Pi$ be the space of probability measure on $\X\times\M\times
  \P{\X\times\M}$. Define $\pi^c_t, \pi^l_t, \pi^g_t$, $t=1,\dots,T$, as
  follows:
  \begin{enumerate}
    \item $\pi^c_t = \PPr(X_t, M_{t-1}, \PXM t{t-1})$.
    \item $\pi^l_t = \PPr(X_t, M_{t-1}, \PXMp t{t-1})$.
    \item $\pi^g_t = \PPr(X_t, M_t, \PXM tt)$.
  \end{enumerate}
\end{definition}
The unconditional PMFs $\pi^c_t, \pi^l_t, \pi^g_t$ defined above are 
information states sufficient for performance evaluation of 
\autoref{prob:structural}. Specifically, they satisfy the following 
properties:
\begin{lemma}
  \label{lemma:info_states}
  $\pi^c_t, \pi^l_t, \pi^g_t$ are information states for the encoder, the 
  memory update and the controller respectively, i.e.,
  \begin{enumerate}
    \item there exist linear transformations $Q^c(c_t)$, $Q^l(l_t)$, and
      $Q^g(g_t)$ such that
      \begin{align}
        \pi^l_t &= Q^c(c_t)\pi^c_t, 
        \label{eq:Qc} \\
        \pi^g_t &= Q^l(l_t)\pi^l_t,
        \label{eq:Ql} \\
        \pi^c_{t+1} &= Q^g(g_t)\pi^g_t.
        \label{eq:Qg}
      \end{align}
    \item the conditional expected instantaneous cost can be expressed as
      \begin{equation}
        \E{\rho(X_t,U_t) | c^t, l^t, g^t} = \tilde\rho(\pi^g_t, g_t),
        \label{eq:expected_cost}
      \end{equation}
      where $\tilde\rho$ is a deterministic function.
  \end{enumerate}
\end{lemma}
\begin{proof}
  This follows from \autoref{lemma:beliefs} and \autoref{def:info_states}.
  See \cite{MahajanTeneketzis:2006b} for details.
  \end{proof}
 Using this result the performance criterion of~\eqref{eq:cost} can 
  be rewritten as
  \begin{align}
    \E{\sum_{t=1}^{T} \rho(X_t, U_t) | \design} &= \sum_{t=1}^{T} \E{ \rho(X_t, 
    U_t) | \design[t]} \nonumber\\
    & \defined \sum_{t=1}^{T}\tilde\rho(\pi^g_t, g_t),
    \label{eq:cost_split}
  \end{align}
  where the sequence $\{\pi^g_1,\dots,\pi^g_T\}$ depends on the choice of 
  (\design). Hence, \autoref{prob:structural} is equivalent to the following 
  deterministic problem:
  \begin{problem}
    \label{prob:deterministic}
    Consider a deterministic system with states $\pi^c_t, \pi^l_t, 
    \pi^g_t$. The initial state $\pi^c_1$ is known and the $t \geq 1$, the 
    system evolves as follows,
    \begin{align}
      \pi^l_t &= Q^c(c_t)\pi^c_t, \label{eq:deter_1}\\
      \pi^g_t &= Q^l(l_t)\pi^l_t, \label{eq:deter_2}\\
      \pi^c_{t+1} &= Q^g(g_t)\pi^g_t\label{eq:deter_3},
    \end{align}
    where $c_t, l_t, g_t$ belong to $\sC, \sL, \sG_s$ respectively and $Q^c, 
    Q^l, Q^g$ are known linear transformations.  At time $t$, an instantaneous 
    cost  $\tilde\rho(\pi^g_t,g_t)$ is incurred.

    The optimization problem is to determine design $(\design)$, where $c 
    \defined (c_1, \dots, c_T)$, $l \defined (l_1, \dots, l_T)$, and $g \defined 
    (g_1, \dots, g_T)$, to minimize the total cost over horizon $T$, i.e.,
    \begin{equation}
      \min_{ (\design) \in \sC^T\times\sL^T\times\sG_s^T} \sum_{t=1}^T \tilde 
      \rho(\pi^g_t, g_t)
    \end{equation}
\end{problem}
This is a classical deterministic optimal control problem; optimal functions 
$(\design[*])$ are determined as follows:
\begin{theorem}
  \label{thm:deterministic}
  An optimal design $(\design[*])$ for \autoref{prob:deterministic} (and 
  consequently for \autoref{prob:structural} and thereby for 
  \autoref{prob:basic}) is given the following nested optimality equations:
  \begin{align}
    V^g_T(\pi^g) &= \inf_{g_T \in \sG_s} \tilde \rho(\pi^g, 
    g_T)\label{eq:DPbasic}\\
    \intertext{and for $t = 1,\dots, T$}
    V^c_t(\pi^c) &= \min_{c_l \in \sC} V^l_t\big( Q^c(c_t) \pi^c \big),
    \label{eq:DPfinite_1} \\
    V^l_t(\pi^l) &= \min_{l_l \in \sL} V^g_t\big( Q^l(l_t) \pi^l \big),
    \label{eq:DPfinite_2}\\
    V^g_t(\pi^g) &= \inf_{g_t \in\sG_s} \tilde \rho(\pi^g,g_t) + V^c_{t+1}\big( Q^g(g_t) \pi^g \big).  
    \label{eq:DPn}
  \end{align}
  The $\arg \min$ (or $\arg \inf$) at each step determines the
  corresponding optimal design for that stage. Furthermore, the optimal
  performance is given by
  \begin{equation}
    \J_T^* = V^c_1(\pi^c_1).
  \end{equation}
\end{theorem}
\begin{proof}
  This is a standard result, see \cite[Chapter~2]{KumarVaraiya:1986}.
\end{proof}

\section{Infinite Horizon Problem}\label{sec:infinite}
In this section we extend the model of \autoref{sec:prob} to an infinite 
horizon ($T \to \infty$) using two performance criteria:
\begin{enumerate}
  \item \emph{Expected Discounted Cost} where the performance of a design is 
    determined by
    \begin{equation}
      \J^\beta(\design) = \E{ \sum_{t=1}^{\infty} \beta^{t-1} \rho(X_t, U_t) 
      | \design },\label{eq:cost_beta}
    \end{equation}
    where $0 < \beta < 1$ is called the discount factor.
  \item \emph{Average Cost per unit time} where the performance of a design 
    is determined by
    \begin{equation}
      \overline \J(\design) = \limsup_{T \to \infty} \frac1T \E{ 
      \sum_{t=1}^{T} \rho(X_t,U_t) | \design}.\label{eq:cost_avg}
    \end{equation}
    We take the $\limsup$ rather than $\lim$ as for some designs $(\design)$ 
    the limit may not exist.
\end{enumerate}
Ideally, while implementing a design for infinite horizon problems, we would 
like to use time-invariant designs. This motivates the following definition.
\begin{definition}
  A design $(\design)$, $c \defined (c_1, c_2, \dots), l \defined (l_1, l_2, 
  \dots), g \defined (g_1, g_2, \dots)$ is called stationary (or 
  time-invariant) if $c_1 = c_2 = \dots = c, l_1 = l_2 = \dots = l, g_1 = g_2 
  = \dots = g$.
\end{definition}

Due to the dynamic team nature of the problem, it is not immediately clear
whether there exist stationary designs that are optimal (or
$\varepsilon$-optimal).  In this section we show that for the expected
discounted cost problem, without loss of optimality, there exist stationary
design that are optimal; for the average cost per unit time problem, under
certain conditions, there exist stationary designs that are
$\varepsilon$-optimal.

\subsection{Expected Discounted Cost Problem} Consider the infinite horizon
problem with  expected discounted cost criterion given
by~\eqref{eq:cost_beta}. For this problem the relations of
\autoref{lemma:beliefs} hold, hence the structural result of
\autoref{thm:structural} is valid, and we can restrict attention to
encoders belonging to $\sG_s$. Define $\pi^c_t, \pi^l_t, \pi^g_t$ as in
\autoref{def:info_states}. \autoref{lemma:info_states} can be proved as
before. The transformations $Q^c, Q^l, Q^g$ and the expected instantaneous
cost $\tilde\rho$ are the same as in the finite horizon case.  Hence, the
infinite horizon problem with the expected discounted cost criterion given
by~\eqref{eq:cost_beta} is equivalent to 
\autoref{prob:deterministic} with the optimization criterion given by
  \begin{equation}
    \J^\beta(\design) \defined \E{ \sum_{t=1}^{\infty} \beta^{t-1}
    \rho(\pi^g_t, g_t) | \design}.
    \label{eq:deterministic_beta_cost}
  \end{equation}

For this problem we have the following result:

\begin{theorem}
  \label{thm:discounted}
  For the infinite horizon expected
  discounted cost problem with the performance criterion given
  by~\eqref{eq:cost_beta}, without loss of optimality, one can restrict
  attention to stationary designs. Specifically, for any optimal design
  (\design[\prime]) there exists a stationary design
  $(c_0^\infty,l_0^\infty,g_0^\infty)$, $c_0^\infty \defined (c_0,c_0,\dots)$,
  $l_0^\infty \defined (l_0,l_0,\dots)$, and  $g_0^\infty \defined
  (g_0,g_0,\dots)$, such that
  \begin{equation}
    V(\pi^c_1) = \J^\beta(c_0^\infty,l_0^\infty,g_0^\infty) = 
    \J^\beta(\design[\prime]),
  \end{equation}
  where $V$ is the unique uniformly bounded fixed point of
  \begin{equation}\label{eq:fixed_point_beta}
    V(\pi) = \min_{(\design) \in \sC\times\sL\times\sG_s}
    \widetilde \rho\big(\widehat Q(c,l)\pi, g\big) +
    \beta V\left( \widetilde Q(\design)(\pi) \right),
  \end{equation}
  with
  \begin{align}
    \widehat Q(c,l) &\defined Q^l(l) \circ Q^c(c) \label{eq:Qhat},\\
    \widetilde Q(\design) &\defined Q^g(g) \circ Q^l(l) \circ Q^c(c),
    \label{eq:Qtilde}
  \end{align}
  and $(c_0,l_0,g_0)$ satisfy
  \begin{equation}
    V(\pi) = \widetilde \rho\left( \widehat Q(c_0,l_0)\pi,g_0 \right) + \beta 
    V\left( \widetilde Q(c_0,l_0,g_0)(\pi) \right).
  \end{equation}
\end{theorem}
\begin{proof}
  See \cite{MahajanTeneketzis:2006b}.
\end{proof}

\subsection{Average Cost per unit time Problem}
Consider the infinite horizon problem with  average cost per unit time 
criterion given by~\eqref{eq:cost_avg}. For this problem the relations of 
\autoref{lemma:beliefs} hold, hence the structural result of 
\autoref{thm:structural} is valid, and we can restrict attention to 
encoders belonging to $\sG_s$. Define $\pi^c_t, \pi^l_t, \pi^g_t$ as in 
\autoref{def:info_states}. \autoref{lemma:info_states} can be proved as 
before. The transformations $Q^c, Q^l, Q^g$ and the expected instantaneous 
cost $\tilde\rho$ are the same as in the finite horizon case.  Hence, the 
infinite horizon expected discounted cost problem is equivalent to 
\autoref{prob:deterministic} with the optimization criterion given by
  \begin{equation}
    \overline{\J}(\design) \defined \limsup_{T \to \infty} \frac 1T\E{ 
    \sum_{t=1}^{T} \rho(\pi^g_t, g_t) | \design}.
    \label{eq:deterministic_avg_cost}
  \end{equation}
For this problem we have the following result:

\begin{theorem}
  \label{thm:average}
  For the infinite horizon average cost
  per unit time problem with the performance criterion given
  by~\eqref{eq:deterministic_avg_cost}, assume
    
  (A1)~for any $\varepsilon > 0$ there exist bounded measurable 
      functions $v(\cdot)$ and $r(\cdot)$ and design $(c_0,l_0,g_0) \in \sC 
      \times \sL \times\sG_s$ such that for all $\pi$,
      \begin{equation}\label{eq:fixed_point_average}
        \begin{split}
          v(\pi) &= \min_{\design \in \sC\times\sL\times\sG_s} v\left( 
          \widetilde Q(\design)\pi \right) \\
          &= v\left( \widetilde Q(c_0,l_0,g_0)\pi \right),
        \end{split}
      \end{equation}
      and
      \begin{multline}
         \min_{\design \in \sC\times\sL\times\sG_s}
        \widetilde \rho\left( \widehat Q(c,l)\pi, g \right) + r\left( 
        \widetilde Q(\design) \pi \right) \leq \\
         \leq v(\pi) + r(\pi)  \leq \widetilde \rho\left( \widehat 
        Q(c_0,l_0)\pi, g_0 \right) + r\left( \widetilde Q(c_0,l_0,g_0)\pi 
        \right) + \varepsilon.
      \end{multline}
  Then for any horizon $T$ and any design (\design[\prime]) for that horizon,
  the  stationary design $(c_0^\infty, l_0^\infty, g_0^\infty)$, $c_0^\infty =
  (c_0, c_0, \dots)$, $l_0^\infty = (l_0,l_0,\dots)$, $g_0^\infty = (g_0, g_0,
  \dots)$, satisfies
  \begin{equation}
    \J_T(c_0^T, l_0^T, g_0^T)  = r(\pi^c_1) + T v(\pi^c_1) \leq 
    \J_T(\design[\prime]) + \varepsilon, \label{eq:diff}
  \end{equation}
  where $\alpha^T = (\alpha, \dots, \alpha)$ ($T$ times) for $\alpha = c_0, 
  l_0, g_0$. Further under (A1), \eqref{eq:diff} is equivalent to
  \begin{equation}
    \overline \J (c_0^\infty, l_0^\infty, g_0^\infty) = v(\pi^c_1) \leq 
    \underline \J (\design[\prime\prime]),\label{eq:avg_inf}
  \end{equation}
  where (\design[\prime\prime]) is any infinite horizon policy and
  \begin{equation}
    \underline \J(\design[\prime\prime]) = \liminf_{T \to \infty} \frac1T 
    \sum_{t=1}^{T} \widetilde \rho\left( \widetilde Q(c''_t,l''_t)\pi_t, g''_t 
    \right).
  \end{equation}
\end{theorem}
\begin{proof}
  See \cite{MahajanTeneketzis:2006b}.
\end{proof}

\subsection{Implication of the Result}
We have shown that there exist optimal stationary designs for the infinite
horizon expected discounted cost problem and, under certain conditions,
$\varepsilon$-optimal stationary designs for the infinite horizon average cost
per unit time problem.  This simplifies the off-line optimization problem
since we have to choose the best amongst $\sC\times\sL\times\sG_s$ stationary
strategies, rather than to choose the best amongst
$\sC^\infty\times\sL^\infty\times\sG_s^\infty$ time-varying strategies.
Further, implementing stationary strategies involves implementing one function
at each agent which is much simpler than implementing a time-varying strategy.

\section{Computational Issues}\label{sec:computation}

The dynamic program of \autoref{thm:deterministic} for joint optimization of encoding, memory updating and
control strategies is similar to a
dynamic program for \pomdp with uncountable state space and uncountable
action space. The information state $\pi^c_t$ belongs to $\Pi$, the space of
probability measures on $\X\times\M\times\P{\X\times\M}$, which is a subset
of probability measures on $\RR^{d}$, with $d = \abs\X\times\abs\M + 1$. The action spaces $\sF$ and $\sL$ are
finite while the action space $\sG_S$ is uncountable. Therefore, in the dynamic
program of \autoref{thm:deterministic}, the information state belongs to the
space of probability measures on a finite dimensional Euclidean space and an
uncountable state space. The standard computational techniques for solving
such \pomdp{s} can be used to obtain numerical results. 

It is the off-line computation of an optimal design that has exponential
complexity. The on-line implementation is simple as we need to implement a stationary design.

\section{Sensors with Imperfect Observations}\label{sec:partial}

So far we have assumed that the sensor perfectly observes the state of the 
system. However in many practical systems, the sensor observations are noisy 
due to external disturbances and the intrinsic noise in the measurement 
hardware. In this section we model this scenario and show that noisy 
observations by the sensor do not alter the nature of the problem. We first 
consider the finite horizon case.

\subsection{Problem Formulation}
\begin{figure}[tb!h]
  \centering
  \resizebox{0.95\linewidth}{!}{
  \empuse{imperfect}}
    \caption{Feedback control system with noisy communication and imperfect 
    observations.}
  \label{fig:imperfect}
\end{figure}
Consider a discrete time imperfect observation system as shown in 
\autoref{fig:imperfect} which operates for $T$ time steps. The state of
the system $X_t$ evolves 
according to~\eqref{eq:state}.  The observations $S_t$ made by the observer 
at time $t$ are noisy version of the state of the system and are given by
\begin{equation}
  S_t = \hat h(X_t, \widehat N_t),
  \label{eq:observations}
\end{equation}
where $\widehat N_t$ denotes the observation noise and $\hat h$ is the 
\emph{observation channel}. $S_t$ takes values in $\S \defined \{1,\dots, 
\abs\S\}$ and $\widehat N_t$ takes values in $\widehat N \defined \{1, \dots, 
\abs{\widehat N}\}$. The sequence of random variables  $\widehat N_1, \dots, 
\widehat N_T$ are i.i.d.\@ with PMF $P_{\widehat N}$.  $\widehat N_1, \dots, 
\widehat N_T$ is also independent of $X_1, W_1, \dots, W_T, N_1, \dots, N_T$.

The sensor is modeled as in \autoref{sec:prob} and operates as follows
\begin{align}
  Z_t &= c_t(S_t,M_{t-1}), \label{eq:obs_enc} \\
  M_t &= l_t(S_t, M_{t-1}). \label{eq:obs_mem}
\end{align}
All other components of the system (the channel, the controller and the 
performance measure) are modeled as in \autoref{sec:prob}. The collection of 
$\IOsystem$ is called an \emph{imperfect observation system}. The 
choice of (\design), $c \defined (c_1,\dots,c_T)$, $l \defined (l_1, \dots, 
l_T)$, $g \defined (g_1, \dots, g_T)$, is called a \emph{design}. The performance of a design, quantified by the expected total cost under 
that design, is given by~\eqref{eq:cost}.  We are interested in the following
optimization problem:

\begin{problem}
  \label{prob:IO} Given an imperfect observation system $\IOsystem$, choose a 
  design (\design[*]) such that
  \begin{equation}
    \J_T(\design[*]) = \J^*_T \defined \min_{\design \in \sC^T \times \sL^T 
    \times \sG^T } \J_T(\design),
  \end{equation}
  where $\sC^T \defined \sC \times \dots \times \sC$ ($T$ times), $\sC$ is the 
  space of functions from $\S\times\M$ to $\Z$, $\sL^T \defined \sL \times 
  \dots \times \sL$ ($T$ times), $\sL$ is the space of functions from 
  $\S\times\M$ to $\M$, $\sG^T \defined \sG_1 \times \dots \times \sG_T$, and 
  $\sG_t$ is the space of functions from $\Y^t\times\U^{t-1}$ to $\U$.
\end{problem}

Although in \autoref{prob:IO} the encoder does not know the state of the
plant, the problem is conceptually same as \autoref{prob:basic} and the
solution methodology of \autoref{prob:basic} works for \autoref{prob:IO} with
very minor changes.

\subsection{Structural Results}

In this section we present structural properties of optimal controllers.  For 
this purpose define the following:
\begin{definition}
  Let \PXM tt and  \PXM {t+1}t be defined as in \autoref{def:info_states}.  
  Define \PXSM {t+1}tt as follows:
  \begin{multline*}
    \PXSM {t+1}tt(x,s,m) \\
    = \Pr{X_{t+1} = x, S_t =s, M_t = m | Y^t, U^t, c^{t+1}, l^t}.
  \end{multline*}
\end{definition}
These beliefs are related as follows:
\begin{lemma}
  For each stage $t$,
    there exists a deterministic functions $\psi$, $\hat \phi$, and
      $\hat \nu$  such that
      \begin{align}
        \PXM {t+1}t &= \psi( \PXM tt, U_t ), \\
        \PXSM {t+1}tt &= \hat\phi( \PXM {t+1}t, c_{t+1} ), \\
        \PXM {t+1}{t+1} &= \hat\nu( \PXSM {t+1}tt, l_{t+1} ).
      \end{align}
\end{lemma}
\begin{proof}
  This can be proved along the same lines as the proof of 
  \autoref{lemma:beliefs}.
\end{proof}
Using the above relationship it can be shown that the structural result of 
\autoref{thm:structural} also hold for \autoref{prob:IO}. Thus, without loss 
of optimality, we can restrict attention to controllers of the 
form~\eqref{eq:structural}. These structural results imply that we can 
formulate a problem equivalent to \autoref{prob:IO} with a time invariant 
action space.

\subsection{Joint Optimization}
We follow the philosophy of \autoref{sec:joint} and use the structural 
results of previous section to obtain a sequential decomposition for 
\autoref{prob:IO}.

\begin{definition}
  Let $\pi^c_t, \pi^l_t, \pi^g_t$, $t=1,\dots,T$ be defined as follows:
  \begin{enumerate}
    \item $\pi^c_t = \Pr{X_t, M_{t-1}, \PXM t{t-1}}$.
    \item $\pi^l_t = \Pr{X_t, M_{t-1}, \PXSM t{t-1}{t-1}}$.
    \item $\pi^g_t = \Pr{X_t, M_t,     \PXM tt}$.
  \end{enumerate}
\end{definition}
\autoref{lemma:info_states} holds for $\pi^c_t, \pi^l_t, \pi^g_t$ defined 
above. Thus, the above unconditional PMFs are information states sufficient
for performance evaluation of  
\autoref{prob:IO}.  This can be shown along the same lines as the proof of 
\autoref{lemma:info_states}. Hence, \autoref{prob:IO} is 
equivalent to a deterministic problem similar to \autoref{prob:deterministic} 
with the transformations $Q^c, Q^l, Q^g$ appropriately defined. The solution 
of this deterministic problem is given by nested optimality equations similar 
to \autoref{thm:deterministic}. Hence, we obtain a sequential decomposition 
of \autoref{prob:IO}. Similar results extend to infinite horizon problems 
using the ideas of \autoref{sec:infinite}.

\section{Conclusion}\label{sec:conclusion}
We have presented a methodology for determining jointly optimal encoding and
control strategies for feedback control systems with limited communication
over noisy channel. The methodology is applicable to finite horizon problems
with expected total cost criterion, to infinite horizon problem with expected
discounted cost criterion, and to infinite horizon problem with average cost
per unit time criterion. We extend this methodology to problem where the
encoder/sensor makes imperfect observations of the state of the system. The
resulting optimality equations can be viewed as \pomdp{s} where the state space
is a real valued vector and the action space is uncountable. Hence traditional
method for solving such \pomdp{s} can be used to obtain a solution for
feedback control problems with communication constraints.

The methodology presented here can be used to obtain a sequential
decomposition of general dynamic team problems with non-classical information
structures.


\bibliographystyle{IEEEtran}
\bibliography{IEEEabrv,../../collection}
\end{document}